\documentclass[12pt]{article}
\usepackage{amsthm, amsmath, amssymb}
\usepackage[version-1-compatibility]{siunitx}
\usepackage{graphicx}
\usepackage{caption}
\usepackage{wrapfig}
\usepackage{tikz}
\usetikzlibrary{decorations.pathreplacing}
\usepackage{rotating}
\usepackage{diagbox}
\usepackage{lscape}
\usepackage{subcaption}
\usepackage{adjustbox}
\usepackage[backend=biber, sorting=none, url=false, style=numeric-comp]{biblatex}
\numberwithin{equation}{section}

\DeclareMathOperator{\Li}{Li}
\DeclareMathOperator{\Ti}{Ti}

\DeclareMathOperator{\csch}{csch}
\DeclareMathOperator{\sech}{sech}
\DeclareMathOperator{\arccot}{arccot}
\DeclareMathOperator{\arccsc}{arccsc}
\DeclareMathOperator{\arcSec}{arcsec}  

\DeclareMathOperator{\arcsinh}{arcsinh}
\DeclareMathOperator{\arccosh}{arccosh}
\DeclareMathOperator{\arctanh}{arctanh}
\DeclareMathOperator{\arccsch}{arccsch}

\DeclareMathOperator{\arccoth}{arccoth}

\graphicspath{ {./} }

\begin{document}

\title{
{Unifying trigonometric and hyperbolic function derivatives via negative integer order polylogarithms} \\
\author{Andrew Ducharme}
\date{February 1, 2024}
}
\maketitle

\begin{centering}
Department of Physics, University of Oregon, 1371 E 13th Avenue, Eugene, Oregon 97403, USA, aducharm@uoregon.edu
\end{centering}

\begin{abstract}
Special functions like the polygamma, Hurwitz zeta, and Lerch zeta functions have sporadically been connected with the $n$th derivatives of trigonometric functions. We show the polylogarithm $\Li_s(z)$, a function of complex argument and order $z$ and $s$, encodes the $n$th derivatives of the cotangent, tangent, cosecant and secant functions, and their hyperbolic equivalents, at negative integer orders $s = -n$. We then show how at the same orders, the polylogarithm represents the $n$th application of the operator $x \frac{d}{dx}$ on the inverse trigonometric and hyperbolic functions. Finally, we construct a sum relating two polylogarithms of order $-n$ to a linear combination of polylogarithms of orders $s = 0,-1,-2, ..., -n$.
\end{abstract}

\section{Introduction}

The polylogarithm, or polylog, is defined as

\begin{equation}
\Li_s(z) = \sum_{k=1}^\infty \frac{z^k}{k^s}.
\end{equation}

\noindent The order and argument $s$ and $z$ can be complex. The function converges absolutely inside the unit disk $|z| \leq 1$ (except at $z = 1$ for $\Re(s) \leq 1$) and can be extended to $|z| > 1$ via analytic continuation. It is principally valued on the complex plane with a branch cut along the positive real axis $1 \leq \Re(z) < \infty$. 

Higher order polylogarithms, like the dilogarithm $\Li_2(z)$ and trilogarithm $\Li_3(z)$, which cannot be written in terms of elementary functions, have been of interest since at least the mid-eighteenth century \cite{maximon_dilogarithm_2003}. Conversely, the negative integer order polylogarithms are simply rational functions in $z$

\begin{equation} \label{negPolyDef}
\Li_{-n} (z) = \left( z \frac{d}{dz} \right)^n \frac{z}{1-z} = \sum_{k=0}^n k! { n+1 \brace k+1} \left(\frac{z}{1-z} \right)^{k+1}
\end{equation}

\noindent where ${n \brace k}$ are Stirling numbers of the second kind. 
Its simplicity likely explains the little discussion of these cases in the literature; Kummer was the only scholar of the polylogarithm Truesdell notes who considered the function at negative integer orders \cite{truesdell_function_1945}. These cases also receive no mention in Lewin's authoritative text on the subject \cite{lewin_dilogarithms_1958}. Negative order polylogarithms were reintroduced by Wood \cite{david_c_wood_computation_1992}, and found by Lee \cite{lee_polylogarithms_1997}, who deemed them ``polypseudologs," in his study of statistical mechanics. Cvijovic derived multiple explicit formulas, one novel, of these functions \cite{cvijovic_polypseudologarithms_2010}. 

Here, we demonstrate how the polylogarithm at negative integer orders is a unifying structure behind the $n$th derivatives of the cotangent, tangent, cosecant, and secant functions, and their hyperbolic cousins. The polylogarithm also underlies the action of the operator $\left( x \frac{d}{dz} \right)^n$ on all twelve of the inverse trigonometric and hyperbolic functions. These results are easy to state: simply check Section \ref{trigDeriv} for the elementary trigonometric derivatives barring $\sin x$ and $\cos x$, Section \ref{hyperDeriv} for the hyperbolic derivatives barring $\sinh x$ and $\cosh x$, and Section \ref{inverseDeriv} for all the inverse trigonometric and hyperbolic functions. Eq. \ref{negPolyDef} simplifies the derivation as well. Simply find the proper substitution, and voil\'a. What sets the polylogarithm apart is its flexibility to represent these 20 functions more simply than other special functions like the Lerch zeta or polygamma functions.

Consider Eq. \ref{negPolyDef} under the transformation $z \rightarrow e^{ix}$, $x \in \mathbb{R}/\{2\pi \ell \}$, $\ell \in \mathbb{Z}$. Then $\frac{d}{dz} \rightarrow (-i)^n \frac{d}{dx}$, $\Li_0(z) = \frac{z}{1-z} \rightarrow \frac{1}{2} \left(i \cot \frac{x}{2} - 1 \right)$, and for $n \geq 1$,

\begin{multline} \label{pl_cot}
\Li_{-n} (e^{ix}) = \frac{i^{1-n}}{2} \left( \frac{d}{dx} \right)^n \cot \frac{x}{2} \\
= \sum_{k=0}^n \frac{k!}{2^{k+1}} {n+1 \brace k+1} \left(i \cot \frac{x}{2} - 1 \right)^{k+1}
\end{multline}

\noindent This is just one of many connections of $n$th cotangent derivative with special functions. Apostol wrote it in terms of the Lerch zeta function $\ell_{1-n}(z) = \frac{i}{2 (2\pi i)^{n-1}} \frac{d^{n-1}}{dz^{n-1}} \cot \pi z$ at negative integer orders \cite{apostol_dirichlet_1970}, a special case of the Lerch transcendent \cite{sousa_lerchs_2023}. It can be written as a polygamma reflection formula $\psi_n (x) - (-1)^n \psi_n (1-x) = -\pi \left( \frac{d}{dx} \right)^n \cot \pi x$ \cite{ks_kolbig_polygamma_1996, boyadzhiev_derivative_2007}, or a Hurwitz zeta reflection formula $\zeta(n, x) + (-1)^n \zeta(n, 1-x)$ \\ $= - \frac{\pi}{(n-1)!} \left( \frac{d}{dx} \right)^{n-1} \cot \pi x$ \cite{adamchik_hurwitz_2007}. But note how Eq. \ref{pl_cot} more succinctly represents the cotangent derivative. Two instances of the polygamma or Hurwitz zeta functions are needed where one polylogarithm would suffice. Only one Lerch transcendent is required, but with more degrees of freedom (three parameters) than the polylogarithm (two parameters). 


The Lerch zeta function and the polylogarithm comparably represent the cotangent derivatives because the polylogarithm substitution $z \rightarrow e^{ix}$ nearly equates the two special functions: $\ell_s(z) = \sum_{k=1}^\infty \frac{e^{2\pi i k z}}{k^{s}} = \Li_s(e^{2\pi i z})$. But the Lerch zeta function is too overconstrained to represent inverse trigonometric and hyperbolic function derivatives. For example, $\Li_{-n}(x) - \Li_{-n}(x) = 2 \left( x \frac{d}{dx} \right)^{n+1} \arctanh x$. To write the left-hand side with Lerch zeta functions, we must undo the complex exponential with a natural logarithm to get a gnarly $\ell_{s} (\frac{1}{2\pi i} \log x )$. But the insertion of the logarithm restricts the Lerch zeta function to exist only over the positive reals, whereas the domain of $\Li_{-n}(x)$ is $x \in \mathbb{R}/\{1\}$. With only one Lerch zeta function per side of the origin, it cannot represent the $\arctanh x$ derivative. The polylogarithm prevails. 

We note the $n$th trigonometric derivatives were connected to negative integer order polylogarithms in \cite{sousa_lerchs_2023} and negative integer order Lerch zeta and Legendre chi functions in \cite{cvijovic_lerch_2010}.

The rest of the paper works through the various trigonometric, hyperbolic, inverse trigonometric, and inverse hyperbolic derivatives. As the difference of mirrored polylogarithms becomes increasingly common in these relations, in Section \ref{chiTi}, we define the Legendre chi function $\chi_s(z)$ and the inverse tangent integral $\Ti_s(z)$, special functions which abbreviate these differences. In Section \ref{pseudoSum}, we use trigonometric derivative polynomials to construct a sum relating two polylogarithms of order $-n$ with a linear combination of all polylogarithms of lesser magnitude negative integer orders.

\section{Trigonometric Derivative Polynomials} \label{trigDeriv}

Trigonometric derivatives and their resulting derivative polynomials have long been a problem of interest \cite{schwatt_introduction_1924, knuth_computation_1967, hoffman_derivative_1995}, and received a flurry of attention in the late 2000s \cite{adamchik_hurwitz_2007, boyadzhiev_derivative_2007, franssens_functions_2007, chang_central_2009, cvijovic_derivative_2009}

To get the $n$th tangent derivative, substitute $z \rightarrow -e^{i x}$ into Eq. \ref{negPolyDef} so for $n \geq 1$,

\begin{equation} \label{pl_tan}
-\Li_{-n}(-e^{ix}) = \frac{i^{1-n}}{2} \left( \frac{d}{dx} \right)^n \tan \frac{x}{2} 
\end{equation}

\begin{equation}
\left( \frac{d}{dx} \right)^n \tan x = 2^n i^{n-1} \sum_{k=0}^n \frac{(-1)^k k!}{2^k} {n+1 \brace k+1} \left(1 + i \tan x \right)^{k+1}.
\end{equation}

The sum of Eqs. \ref{pl_cot} and \ref{pl_tan} is proportional to the $n$th derivative of $\cot \frac{x}{2} + \tan \frac{x}{2} = 2\csc x$, so

\begin{equation} \label{csc}
\Li_{-n}(e^{ix}) - \Li_{-n}(-e^{ix}) = i^{1-n} \left( \frac{d}{dx} \right)^n \csc x.
\end{equation}

\noindent Expanding both $(1 + i \tan x)^{k+1}$ and $(i \cot x - 1)^{k+1}$ by the Binomial Theorem gives

\begin{multline} \label{Dcsc}
\left( \frac{d}{dx} \right)^n \csc x = \frac{i^{n-1}}{2} \sum_{k=0}^n \frac{(-1)^k k!}{2^k} {n+1 \brace k+1} \\
\times \sum_{j=0}^{k+1} {k+1 \choose j} i^j \left[\tan^j \frac{x}{2} - (-1)^j \cot^j \frac{x}{2} \right].
\end{multline}

\noindent This is similar to the polylogarithm duplication formula $\Li_s(z) + \Li_s(-z) = 2^{1-s} \Li_s(z^2)$ which under the same substitution $z \rightarrow e^{ix}$ proves a cotangent double angle formula $2 \cot 2x = \cot x - \tan x$.

Taking advantage of $\csc \left( x + \frac{\pi}{2} \right) = \sec x$, we get our final trig derivative

\begin{equation}
\Li_{-n}(ie^{ix}) - \Li_{-n}(-ie^{ix}) = i^{1-n} \left( \frac{d}{dx} \right)^n \sec x.
\end{equation}

\noindent Because the identities for $\tan \left( \theta + \frac{\pi}{4} \right) = \frac{\tan \theta + 1}{1 - \tan \theta}$ and $\cot \left( \theta + \frac{\pi}{4} \right) = \frac{\cot \theta - 1}{1 + \cot \theta}$ are not simply linear combinations of other trig functions, the $n$th secant derivative can be written more succinctly as Eq. \ref{Dcsc} with $\frac{x}{2} \rightarrow \frac{x}{2} + \frac{\pi}{4}$. Alternatively, combining the fractional forms of the quarter-period-shifted tangent and cotangent, then binomially expanding $(\sec x \pm \tan x)^j$, produces the following three sum equation

\begin{multline} \label{Dsec}
2i^{1-n} \left( \frac{d}{dx} \right)^n \sec x = \\ 
\sum_{k=0}^n \frac{(-1)^k k!}{2^k} {n+1 \brace k+1} \sum_{j=0}^{k+1}  i^j {k+1 \choose j} \sum_{\ell = 0}^j {j \choose \ell} \tan^{j-\ell}x \sec^\ell x (1 - (-1)^\ell).
\end{multline}

Two polylogarithms is likely the fewest needed to build the cosecant and secant derivatives because one, these even functions are constructed from the odd cotangent and tangent, and two, there is no generic identity to reduce the form $\Li_s(z) - \Li_s(-z)$, as it itself is a special function, the Legendre chi function $\chi_s(x)$.

\section{Legendre chi function and Inverse tangent integral} \label{chiTi}

As infinite sums, the Legendre chi function

\begin{equation}
\chi_s(z) = \sum_{k=0}^\infty \frac{z^{2k+1}}{(2k+1)^s} = \frac{1}{2} [ \Li_s(z) - \Li_s(-z) ]
\end{equation}

\noindent and inverse tangent integral

\begin{equation}
\Ti_s(z) = \sum_{k=0}^\infty (-1)^k \frac{z^{2k+1}}{(2k+1)^s} = \frac{1}{2i} [ \Li_s(iz) - \Li_s(-iz) ]
\end{equation}

\noindent converge on the unit disk for all $s$, except at $z = \pm1$ for $\chi_s(z)$, and $z = \pm i$ for $\Ti_s(z)$ when $s \leq 1$. The functions are the real and imaginary parts of the polylogarithm, and are related by $\Ti_s(z) = -i \chi_s(iz)$. The inverse tangent integral is so named because $\Ti_2(z) = \int_0^x \frac{\arctan t}{t} dt$. Consequently, $\Ti_1(z) = \arctan z$ and $\chi_1(z) = \arctanh z$. The first few functions of zeroth and negative integer order are 

\begin{table}[h]
\begin{adjustbox}{width = 0.75\textwidth, center = \textwidth}
\begin{tabular}{c|cc}
$n$ & $\chi_{n}(z)$ & $\Ti_n(z)$ \\ \hline
0 & $\frac{z}{1-z^2}$ & $\frac{z}{1 + z^2}$ \\
-1 & $\frac{z + z^3}{(1-z^2)^2}$ & $\frac{z - z^3}{(1 + z^2)^2}$ \\
-2 & $\frac{z + 6z^3 + z^5}{(1-z^2)^3}$ & $\frac{z - 6z^3 + z^5}{(1 + z^2)^3}$ \\
-3 & $\frac{z + 23z^3 + 23z^5 + z^7}{(1-z^2)^4}$ & $\frac{z - 23z^3 + 23z^5 - z^7}{(1 + z^2)^4}$ \\
-4 & $\frac{z + 76z^3 + 230 z^5 + 76 z^7 + z^9}{(1-z^2)^5}$  & $\frac{z - 76z^3 + 230 z^5 - 76 z^7 + z^9}{(1 + z^2)^5}$
\end{tabular}
\end{adjustbox}
\end{table}

A check of the OEIS finds that the coefficients of the numerator polynomial in the Legendre chi function $S(n, k)$ is sequence A060187, describing the Eulerian numbers of type B \cite{oeis}. This pairs nicely with the similar coefficients in the numerator of the polylogarithm being the Eulerian numbers of type A (often just called the Eulerian numbers). A single sum form provided by the OEIS is 

\begin{equation}
S(n,k) = \sum_{j = 1}^k (-1)^{k-j} {n + 1 \choose k-j} (2j-1)^{n}.
\end{equation}

Closed forms of the Legendre chi function and inverse tangent integral over non-positive integers $n \leq 0$ are

\begin{equation}
\chi_{-n} (z) = \frac{1}{(1-z^2)^{n+1}} \sum_{k=1}^{n+1} S(n,k) z^{2k-1}
\end{equation}

\begin{equation}
\Ti_{-n} (z) = -\frac{1}{(1+z^2)^{n+1}} \sum_{k=1}^{n+1} (-1)^k S (n,k) z^{2k-1}.
\end{equation}

Using these to clean up Eqs. \ref{Dcsc} and \ref{Dsec}, we get

\begin{equation}
\left( \frac{d}{dx} \right)^n \csc x = \frac{(-1)^n}{2^n} e^{-2ix} \csc^{n+1} x  \sum_{k=1}^{n+1} S(n,k) e^{-i(n-2k)x}
\end{equation}

\begin{equation}
\left( \frac{d}{dx} \right)^n \sec x = -\frac{i^n}{2^n} e^{-2ix} \sec^{n+1} x  \sum_{k=1}^{n+1} (-1)^k S(n,k) e^{-i(n-2k)x}
\end{equation}

\section{Hyperbolic Function Derivatives} \label{hyperDeriv}

The hyperbolic functions have similar, but simpler, derivatives compared to the trigonometric functions, lacking coefficients with $n$ dependence because the substitution of the real exponential takes $\left( z \frac{d}{dz} \right)^n$ to $\left( e^x e^{-x} \frac{d}{dx} \right)^n = \frac{d^n}{dx^n}$, with no extra factor of $-i^n$. These equations hold over $x \in \mathbb{R}$.

\begin{equation}
\Li_{-n}(e^{x}) = - \frac{1}{2} \left( \frac{d}{dx} \right)^n \coth \frac{x}{2}
\end{equation}

\begin{equation}
\Li_{-n}(-e^{x}) = - \frac{1}{2} \left( \frac{d}{dx} \right)^n \tanh \frac{x}{2}
\end{equation}

\begin{equation}
\Li_{-n}(e^{x}) - \Li_{-n}(-e^{x}) = 2 \chi_{-n}(e^x) =  - \left( \frac{d}{dx} \right)^n \csch x,
\end{equation}

\begin{equation}
\Li_{-n}(ie^{x}) - \Li_{-n}(-ie^{x}) = 2i \Ti_{-n}(e^x) = i \left( \frac{d}{dx} \right)^n \sech x,
\end{equation}

\begin{equation}
\left( \frac{d}{dx} \right)^n \coth x
= 2^n \sum_{k=0}^n \frac{(-1)^k k!}{2^{k}} {n+1 \brace k+1} \left( 1 + \coth x \right)^{k + 1}
\end{equation}

\begin{equation}
\left( \frac{d}{dx} \right)^n \tanh x
= 2^n \sum_{k=0}^n \frac{(-1)^k k!}{2^{k}} {n+1 \brace k+1} \left( 1 + \tanh x \right)^{k + 1}
\end{equation}

\begin{equation}
\left( \frac{d}{dx} \right)^n \csch x = \frac{(-1)^n}{2^n} e^{2x} \csch^{n+1} x  \sum_{k=1}^{n+1} S(n,k) e^{(n-2k)x}
\end{equation}

\begin{equation}
\left( \frac{d}{dx} \right)^n \sech x = -\frac{(-1)^n}{2^n} e^{2x} \sech^{n+1} x  \sum_{k=1}^{n+1} (-1)^k S(n,k) e^{(n-2k)x}
\end{equation}

\section{Inverse Trig and Hyperbolic Functions} \label{inverseDeriv}

The substitution of exponentials in the previous sections are useful because it reduces the operator $z \frac{d}{dz}$ to a derivative, but are not the sole use cases of Eq. \ref{negPolyDef}. We can build the $n$th application of $x \frac{d}{dx}$ to inverse trigonometric and hyperbolic functions from $\Ti_0 (x) = x \frac{d}{dx} \arctan x$ and $\chi_0(x) = x \frac{d}{dx} \arctanh x$ for $n \geq 0$ and, unless otherwise noted, $x \in \mathbb{R}$.

\begin{equation}
\chi_{-n}(x) = \left( x \frac{d}{dx} \right)^{n+1} \arctanh x
\end{equation}

\begin{equation}
\chi_{-n} \left( \frac{1}{x} \right) = \left( -x \frac{d}{dx} \right)^{n+1} \arccoth x
\end{equation}

\begin{equation}
\chi_{-n} \left( \frac{x}{\sqrt{1+x^2}} \right) = \left( (x + x^3) \frac{d}{dx} \right)^{n+1} \arcsinh x
\end{equation}

\begin{equation}
\chi_{-n} \left( \frac{1}{x \sqrt{1+x^{-2}}} \right) = \left( -\left( x + \frac{1}{x} \right) \frac{d}{dx} \right)^{n+1} \arccsch x
\end{equation}

\noindent Note for $x < 0$, $x\sqrt{1+x^{-2}} = -\sqrt{1+x^2}$.

\begin{equation}
\chi_{-n} \left( \frac{\sqrt{x^2 - 1}}{x} \right) = \left( x \left( x^2 - 1 \right) \frac{d}{dx} \right)^{n+1} \arccosh x, \; x \geq 1
\end{equation}

\begin{equation}
\chi_{-n} \left( \sqrt{1 - x^2} \right) = \left( \left( x - \frac{1}{x} \right) \frac{d}{dx} \right)^{n+1} \arccosh x, \; 0 < x \leq 1.
\end{equation}

\begin{equation}
\Ti_{-n}(x) = \left( x \frac{d}{dx} \right)^{n+1} \arctan x 
\end{equation}

\begin{equation}
\Ti_{-n}(x) = - \left( x \frac{d}{dx} \right)^{n+1} \arccot x
\end{equation}

\begin{equation}
\Ti_{-n} \left( \frac{x}{\sqrt{1-x^2}} \right) = \left( (x - x^3) \frac{d}{dx} \right)^{n+1} \arcsin x, \; |x| \leq 1
\end{equation}

\begin{equation}
\Ti_{-n} \left( \frac{\sqrt{1-x^2}}{x} \right) = \left( (x^3 - x) \frac{d}{dx} \right)^{n+1} \arccos x, \; |x| \leq 1
\end{equation}

\begin{equation}
\Ti_{-n} \left( \frac{1}{x \sqrt{1-x^{-2}}} \right) = \left( \left( \frac{1}{x} - x \right) \frac{d}{dx} \right)^{n+1} \arccsc x, \; |x| \geq 1
\end{equation}

\begin{equation}
\Ti_{-n} \left( x \sqrt{1-x^{-2}}\right) = \left( \left( x - \frac{1}{x} \right) \frac{d}{dx} \right)^{n+1} \arcSec x, \; |x| \geq 1
\end{equation}

Alternatively, one could simplify the operand of $\left( z \frac{d}{dz} \right)^n$ at the expense of complicating the operator. To transform the operand into a generic $f(x)$, choose $z \rightarrow \frac{f(x)}{1+f(x)}$ so

\begin{multline}
 \Li_{-n} \left( \frac{f(x)}{1+f(x)} \right) = \left( \frac{f(x)}{f'(x)} (1+f(x)) \frac{d}{dx} \right)^n f(x) \\
  = \sum_{k=0}^n k! {n+1 \brace k+1} f(x)^{k+1}.
\end{multline}

\noindent $\Li_n(z)$ diverges for $z = 1$, but $\frac{f(x)}{1 + f(x)}$ can never be unitary. $\Li_{-n} \left( \frac{f(x)}{1+f(x)} \right)$ does diverge when $f(x) = -1$.

As an example, we will choose operands $f(x) = \sin x$ and $\cos x$.

\begin{multline}
 \Li_{-n} \left( \frac{\sin x}{1 + \sin x} \right) = \left( \tan x (1 + \sin x) \frac{d}{dx} \right)^n \sin x \\
  = \sum_{k=0}^n k! {n+1 \brace k+1} \sin^{k+1} x
\end{multline}

\begin{multline}
 \Li_{-n} \left( \frac{\cos x}{1 + \cos x} \right) = \left( - \cot x (1 + \cos x) \frac{d}{dx} \right)^n \cos x \\
  = \sum_{k=0}^n k! {n+1 \brace k+1} \cos^{k+1} x.
\end{multline}

\section{A New Ladder-like Sum} \label{pseudoSum}

Polylogarithmic ladders are extraordinary numerical connections between linear combinations of polylogarithms of fixed order $n$ evaluated at different powers of some $x$ with the natural logarithm $\log (y) = -\Li_1(1-y)$. \cite{MR1415794} provides the example

\begin{multline*}
- \pi^2 \log(2) = \frac{35}{2} \Li_3(1) + 36 \Li_3 \left( \frac{1}{2} \right) \\
- 18 \Li_3 \left( \frac{1}{4} \right) - 4 \Li_3 \left( \frac{1}{8} \right) + \Li_3 \left( \frac{1}{64} \right).
\end{multline*}

\noindent We will instead bridge two (fixed) negative integer order polylogarithms $-n$ to a linear combination of individual polylogarithms of each integer order from $s = 0$ to $s = -n$. An illustrating example is

\begin{multline}
\frac{z}{2} \left( \Li_{-6}(z) - \Li_{-6}(-z) \right)= \\
\Li_{0}(z^2) - 12 \Li_{-1}(z^2) + 60 \Li_{-2}(z^2) -  160 \Li_{-3}(z^2)\\
+ 240 \Li_{-4}(z^2) -192 \Li_{-5}(z^2) + 64 \Li_{-6}(z^2).
\end{multline}

We derive the general form by combining Eq. \ref{csc} with an alternate way Risomar Sousa used get the $n$th derivative of $\csc x$ \cite{sousa_lerchs_2023}: apply the General Leibniz rule to $\csc x = e^{-ix} (i + \cot x)$ so

\begin{equation}
\left( \frac{d}{dx} \right)^n \csc x = 2 i^{n-1} e^{-ix} \sum_{k=0}^n {n \choose k} (-1)^{n-k} 2^k \Li_{-n} (e^{2ix}).
\end{equation}

\noindent The combination is

\begin{equation} \label{pseudo}
\Li_{-n}(e^{ix}) - \Li_{-n}(-e^{ix}) = 2 e^{-ix} \sum_{k=0}^n {n \choose k} (-1)^{n-k} 2^k \Li_{-k} (e^{2ix}),
\end{equation}

\noindent which holds for $n \geq 0$. $e^{ix} \rightarrow z$ gives the final relation

\begin{equation} \label{Pseudo}
\Li_{-n}(z) - \Li_{-n}(-z) = \frac{2}{z} \sum_{k=0}^n {n \choose k} (-1)^{n-k} 2^k \Li_{-k} (z^2).
\end{equation}

\noindent This could be interpreted as the negative integer version of 

\begin{multline}
\Li_n (z) = \frac{(-1)^{n-1}}{(n-2)!} \int_1^z \frac{\log (t^{n-2}) \log (1 - t)}{t} dt + \Li_n (1) \\
+ \sum_{k=1}^{n-2} \frac{(-1)^{k-1}}{k!} \Li_{n-k}(z) \log^k (z)
\end{multline}

\noindent from \cite{mathematica}.

The first few sums are

\begin{equation}
\Li_{0}(z) - \Li_{0}(-z) = \frac{2}{z} \Li_0(z^2)
\end{equation}

\begin{equation}
\Li_{-1}(z) - \Li_{-1}(-z) = \frac{2}{z} \left[ -\Li_0(z^2) + 2 \Li_{-1}(z^2) \right]
\end{equation}

\begin{equation}
\Li_{-2}(z) - \Li_{-2}(-z) = \frac{2}{z} \left[ \Li_0(z^2) - 4 \Li_{-1}(z^2) + 4 \Li_{-2}(z^2) \right]
\end{equation}

\begin{multline}
\Li_{-3}(z) - \Li_{-3}(-z) = \\
\frac{2}{z}  \left[-\Li_{0}(z^2) + 6 \Li_{-1}(z^2) - 12 \Li_{-2}(z^2) + 8 \Li_{-3}(z^2) \right]
\end{multline}

\begin{multline}
\Li_{-4}(z) - \Li_{-4}(-z) = \\
\frac{2}{z} \left[\Li_{0}(z^2) - 8 \Li_{-1}(z^2) + 24 \Li_{-2}(z^2) - 32 \Li_{-3}(z^2) + 16 \Li_{-4}(z^2) \right]
\end{multline}

The same approach replicated on $\sec x = e^{-ix}(1 + i \tan x)$ gives a near-identical result

\begin{equation}
\Li_{-n}(ie^{ix}) - \Li_{-n}(-ie^{ix}) = -2 i e^{-ix} \sum_{k=0}^n {n \choose k} (-1)^{n-k} 2^k \Li_{-k} (-e^{2ix}),
\end{equation}

\noindent which is Eq. \ref{pseudo} under $e^{ix} \rightarrow i e^{ix}$ ($z \rightarrow iz$). Writing Eq. \ref{Pseudo} as a Legendre chi function 

\begin{equation}
z \chi_{-n} (z) = \sum_{k=0}^n {n \choose k} (-1)^{n-k} 2^k \Li_{-k} (z^2)
\end{equation}

\noindent demonstrates we have overcomplicated the same identity derivation for the inverse tangent integral

\begin{equation}
z \Ti_{-n} (z) = - \sum_{k=0}^n {n \choose k} (-1)^{n-k} 2^k \Li_{-k} (-z^2).
\end{equation}

\printbibliography

\end{document}